\documentclass{article}
\usepackage{amsmath, amssymb}
\usepackage{latexsym}
\usepackage[utf8]{inputenc}
\usepackage[T1]{fontenc}
\usepackage{geometry}
\usepackage{makeidx}
\usepackage{showkeys}
\makeindex

\newcommand{\R}{{\mathbf R}}      
\newcommand{\N}{{\mathbf N}}
\newcommand{\C}{{\mathbf C}}
\newcommand{\Z}{{\mathbf Z}}
\newcommand{\Q}{{\mathbf Q}}

\newcommand{\bx}{\hfill{$\Box $}}

\newcommand{\Img}{\textrm{Im}\,}

\newcommand{\charac}{\textrm{char}}

\newtheorem{lemma}{Lemma}
\newtheorem{corollary}{Corollary}
\newtheorem{proposition}{Proposition}
\newtheorem{theorem}{Theorem}
\newtheorem{exercise}{Exercise}
\newtheorem{result}{Result}
\title{Hilbertian fields and Hilbert's irreducibility theorem}
\author{Rodney Coleman, Laurent Zwald}
\begin{document}
\maketitle

\begin{abstract}
\noindent Hilbert's irreducibility theorem plays an important role in inverse Galois theory. In this article we introduce Hilbertian fields and present a clear detailed proof of Hilbert's irreducibility theorem in the context of these fields. 
\end{abstract}

An important result in inverse Galois theory is Hilbert's irreducibility theorem. Unfortunately it is difficult to find a clear proof, probably because such a proof requires many detailed steps. Those that we have seen lack important details or have errors and so make reading difficult. For this reason we set out to write a clear, detailed proof, which a reader with a certain mathematical maturity should not find difficult. We will use some basic results from Galois theory, which we detail in an appendix. (The word Result in the text refers to these results.)\\

We begin by introducing Hilbertian fields. Let $f(X,Y)$ be a nonzero polynomial in two variables over a field $F$. Collecting monomials having the same power of $Y$, we may write
$$
f(X,Y) = a_0(X) + a_1(X)Y + a_2(X)Y^2 + \cdots + a_n(X)Y^n,
$$
where the $a_i(X)$ are polynomials in $X$ alone and $a_n(X)\neq 0$, i.e., we may consider $f$ as a member of $F[X][Y]$. The number $n$ is the degree of $f$ with respect to $Y$. We recall that a nonzero element $a$ in an integral domain $R$ is {\it irreducible}\index{irreducible element} if it is a nonunit and, whenever $a=bc$, either $b$ or $c$ is a unit. As a polynomial ring over an integral domain is an integral domain, $F[X][Y]$ is an integral domain. For reasons which will become obvious further on, we will say that $f\in F[X,Y]$ is irreducible, if $f$ is irreducible as an irreducible element of the ring $F[X][Y]$ and has degree greater than $0$ in $Y$, i.e., the polynomial has at least one monomial containing a power of $Y$. If a polynomial is not irreducible, then we will say it is reducible. 

We may extend this definition to polynomials in more than two variables. If $f\in F[X_1, \ldots ,X_k]$, with $k\geq 3$, then we may consider $f$ as an element of $F[X_1,\ldots ,X_{k-1}][X_k]$. We will say that $f$ is irreducible if $f$ is irreducible in the polynomial ring $F[X_1,\ldots ,X_{k-1}][X_k]$ and has degree greater than $0$ in $X_k$. (This definition is not entirely satisfactory, because it depends on which variable we set in the last position.)

Let $f\in F[X,Y]$ be a polynomial of degree greater than $0$ in $Y$. For every $b\in F$, we may define a polynomial $f_b\in F[Y]$ by setting $f_b(Y)=f(b,Y)$. If $a_n(b)\neq 0$, then $f_b$ has $n$ roots, counted according to their multiplicity. If these roots are distincts, then we say that $b$ is a regular value. 

\begin{proposition}\label{propHILBERTreduc1} Let $f(X,Y)$ be a polynomial of degree greater than $0$ in $Y$ over a field $F$ of characteristic $0$. Then all but a finite number of values $b\in F$ are regular.
\end{proposition} 

\noindent \textsc{proof} If we eliminate those values of $b$, which are roots of the leading coefficient $a_n$, then the polynomial $f_b$ has a positive degree. We now consider $f$ as an element of $F(X)[Y]$ and it is not difficult to see that $\Delta (f)(b)=\Delta (f_b)$, where $\Delta (f)$ (resp. $\Delta (f_b)$) denotes the discriminant of $f$ (resp. $f_b$). As $\Delta (f)$ is a polynomial with coefficients in $F(X)$, there is a finite number of elements $u$ of $F(X)$ for which $\Delta (f)(u)=0$; in particular, there is a finite number of values of $b\in F$ for which $\Delta (f)(b)=0$. If exclude these values, then $\Delta (f_b)=\Delta (f)(b)\neq 0$, i.e., $b$ is regular. \bx\\

We may now define the notion of a Hilbertian field. If for any $f\in F[X,Y]$ which is irreducible, there exists an infinite number of values of $b\in F$ such that $f_b(Y)=f(b,Y)\in F[Y]$ is irreducible, then we say that the field $F$ is Hilbertian. We say that $f_b$ is a specialization of $F$. Clearly a finite field cannot be Hilbertian. This is also the case for a field which is algebraically closed. An important example of a Hilbertian field is that of the rational numbers $\Q$. This is known as Hilbert's irreducibility theorem. We will prove this in Section 4. of the article. For the moment we will consider certain important properties of Hilbertian fields.\\ 

\noindent \underline{\bf 1. Properties of Hilbertian fields} \\

In this section we present some technical results, which enable us to illustrate how the notion of a Hilbertian field intervenes in inverse Galois theory. In particular, the aim is to arrive at an important result at the end of the section, namely Theorem \ref{thmHILBERTprop1}. It could be useful to look at this theorem before reading the section in detail, this in order to appreciate the direction of the section. 

\begin{lemma}\label{lemHILBERTprop1} Let $R$ be an integral domain, $S$ a subring of $R$ and $f,h\in S[X]$, with $f$ monic. If $g\in R[X]$ and $fg=h$, then $g\in S[X]$.
\end{lemma}

\noindent \textsc{proof} As $f$ is monic, there exist $q,r\in S[X]$ such that $h=fq+r$, with $\deg r<\deg f$. We have
$$
fq + r =fg \Longrightarrow r=f(g-q).
$$
As $R$ is an integral domain,
$$
\deg r = \deg f + \deg (g-q) \Longrightarrow g-q = 0,
$$
because $\deg r <\deg f$. Therefore $g=q\in S[X]$.\bx\\

The next preliminary result concerns Galois extensions of fields of fractions and is interesting in its own right.

\begin{proposition} \label{proposHILBERTprop1} Suppose that $R$ is an integral domain and $F$ its field of fractions. In addition, let $E$ be a separable extension of $F$ of degree $n$. Then there exists $\alpha \in E$ such that $E=F(\alpha )$ and $m(\alpha ,F)\in R[X]$.
\end{proposition} 

\noindent \textsc{proof} From the primitive element theorem we know that there exists $\beta\in E$ such that $E=F(\beta )$. As $F$ is the field of fractions of $R$, we may multiply $m=m(\beta ,F)$ by a nonzero constant $d\in R$ to obtain $dm\in R[X]$. Setting $\alpha = d\beta $, we have $F(\alpha )=F(\beta )$. We now look for $m(\alpha ,F)$. If
$$
f(X) = d^nb_0 + d^{n-1}b_1X + \cdots + db_{n-1}X^{n-1} + X^n,
$$
where 
$$
m(X) = b_0 + b_1X + \cdots + b_{n-1}X^{n-1} + X^n,
$$
then $f\in R[X]$ and
$$
f(\alpha ) = f(d\beta ) = d^nm(\beta )=0.
$$
Also, $f$ is monic and
$$
[E:F] = [F(\alpha ):F] = [F(\beta ):F] =n
$$
and so $f$ is the minimal polynomial of $\alpha$ over $F$ and $f\in R[X]$. This proves the result.\bx\\

At this point we introduce some notation. If $R$ is a subring of the field $F$ and $A$ a subset of an extension $E$ of $F$, then we will write $R[A]$ for the subring of $E$ generated by $R$ and $A$. If $A$ is composed of a single element $a$, then we will write $R[a]$ for $R[\{a\}]$. From now on we will suppose that rings and fields have characteristic $0$. The next result is fundamental.

\begin{lemma}\label{lemHILBERTprop2} We take $R$, $F$, $E$, $\alpha$ as in Proposition \ref{proposHILBERTprop1}, with $f=m(\alpha ,F)\in R[X]$, and $A$ a finite subset of $E$ containing $\alpha$ such that
$$
\forall x\in A\; \forall \sigma \in Gal(E/F),\; \sigma (x) \in A.
$$
Then there exists $u\in R$ such that, for any field $F'$ and ring homomorphism $\omega :R\longrightarrow F'$, with $\omega (u)\neq 0$, we may find a Galois extension $E'$ of $F'$ and a ring homomorphism extension $\tilde{\omega}:R[A]\longrightarrow E'$ of $\omega$ with the following properties:
\begin{itemize}
\item $E'=F'(\alpha ' )$, where $\alpha '=\tilde{\omega }(\alpha )$;
\item If $f'\in F'[X]$ is the polynomial obtained from $f$ by applying $\omega$ to the coefficients of $f$ and $f'$ is irreducible, then $G'=Gal(E'/F')$ is isomorphic to $G=Gal(E/F)$.
\end{itemize}
\end{lemma}

\noindent \textsc{proof} The proof of this result is rather long, so we will proceed by steps.\\

\noindent \underline{1. Definition of $u$}: Let $u=\Delta (f)$, the discriminant of $f$. (For a definition of the discriminant, see for example \cite{rotman}). As $\charac F=0$, because $\charac R =0$, and $f$ is irreducible, $f$ has no multiple root. This implies that $u\neq 0$. If $F'$ is a field and $\omega :R\longrightarrow F'$ a ring homomorphism such that $\omega (u)\neq 0$, then $\Delta (f') = \omega (u)\neq 0$, hence $f'$ is strongly separable.\\

\noindent \underline{2. A first extension of $R$ and $\omega$}:  We now construct a ring $\tilde{R}$, containing $R$ and we extend $\omega$ to this ring. As $E=F(\alpha )$ and $A\subset E$, for every $x\in A$, there exists $g_x\in F[X]$ such that $x=g_x(\alpha )$. In addition, $F$ is the field of fractions of $R$, and so there exists $d_x\in R^*$ such that $d_xg_x\in R[X]$. We now set
$$
d = \prod _{x\in A} d_x.
$$
Then $dg_x\in R[X]$, for all $x\in A$. We now set
$$
\tilde{R} = R[d^{-1}]\subset F
$$
and extend $\omega$ to $\omega _1:\tilde{R}\longrightarrow F'$ by setting $\omega _1(d^{-1})=\omega (u)^{-1}$.

It should be noticed that $\tilde{R}[A]=\tilde{R}[\alpha]$. First, $\alpha \in A$ implies that $\tilde{R}[\alpha ]\subset \tilde{R}[A ]$. On the other hand, if $x\in A$ and $g_x(X)=\sum _{i=0}^na_iX^i$, with $a_0, \ldots ,a_n\in F$, then 
$$
x = g_x(\alpha ) = \sum _{i=0}^na_i\alpha ^i \Longrightarrow dx = \sum _{i=0}^n\left(\prod _{y\in A, y\neq x}d_y\right)(d_xa_i)\alpha ^i,
$$
which lies in $R[\alpha ]$, because $d_xa_i\in R$, for all $i$, and $d_y\in R$, for all $y$. However, $R[\alpha ]\subset \tilde{R}[\alpha ]$ and 
$x=d^{-1}(dx)\in \tilde{R}[\alpha ]$. Hence, $A\subset \tilde{R}[\alpha ]$ and so $\tilde{R}[A]\subset \tilde{R}[\alpha ]$.\\

\noindent \underline{3. $\tilde{R}[X]/(f)$ and $\tilde{R}[\alpha ]$ are isomorphic}: There is a natural homomorphism from $\tilde{R}[X]$ into $\tilde{R}[\alpha ]$:
$$
\phi :\tilde{R}[X]\longrightarrow \tilde{R}[\alpha ], g \longmapsto g(\alpha ).
$$
If $h\in \ker \phi$, then there exists $g\in F[X]$ such that $h=fg$, because $f=m(\alpha ,F)$. From Lemma \ref{lemHILBERTprop1}, $g\in \tilde{R}[X]$, because $\tilde{R}$ is a subring of $F$. Therefore $\ker \phi \subset (f)$. On the other hand, if $g\in (f)$, then $g(\alpha )=0$ and so $g\in \ker \phi$. It follows that $\ker \phi =(f)$. As $\phi$ is surjective, we have an isomorphism
$$
\bar{\phi} : \frac{\tilde{R}[X]}{(f)} \longrightarrow \tilde{R}[\alpha ].
$$
\noindent \underline{4. Construction of the extension $E'$ of $F'$}: Our next task is to construct a Galois extension $E'$ of $F'$ and a ring homomorphism $\tilde{\omega}$ from $\tilde{R}[A]$ into $E'$, extending $\omega _1$ and hence $\omega$. Let $g'$ be an irreductible factor of $f'$ and $\rho : F'[X]:\longrightarrow F'[X]/(g')$ the natural projection. From the homomorphism $\omega _1:\tilde{R}\longrightarrow F'$ constructed above, we obtain the natural homomorphism $\hat{\omega}_1:\tilde{R}[X]\longrightarrow F'[X]$. We now compose $\hat{\omega}_1$ with $\rho$ to obtain the homomorphism
$$
\rho \circ \hat{\omega}_1: \tilde{R}[X] \longrightarrow F'[X]/(g')
$$
and then use this to define another homomorphism:
$$
\gamma : \frac{\tilde{R}[X]}{(f)} \longrightarrow \frac{F'[X]}{(g')}, v + (f) \longmapsto \rho \circ \hat{\omega}_1(v).
$$
(As $\rho \circ \hat{\omega}_1(f) = f'+(g')$ and $g'\vert f'$, we must have $\rho \circ \hat{\omega}_1(v)=0$, for all $v\in (f)$, hence $\gamma$ is well-defined.)

Now we set 
$$
E' = \frac{F'[X]}{(g')} \qquad \textrm{and} \qquad\tilde{\omega} = \gamma \circ \bar{\phi}^{-1}.
$$
As $g'$ is irreducible $E'$ is a field, which is clearly an extension of $F'$. Also, $\tilde{R}[A]=\tilde{R}[\alpha ]$ and so $\tilde{\omega}$ is a homomorphism from $\tilde{R}[A]$ into $E'$. We need to check that $\tilde{\omega}$ extends $\omega$. If $x\in R\subset \tilde{R}[A]$, then 
\begin{eqnarray*}
\tilde{\omega}(x) &=& \gamma \circ \bar{\phi}^{-1}(x)\;\;=\;\;\gamma (x+(f))\;\;=\;\;\rho \circ \hat{\omega}_1(x)\\
				&=& \rho (\omega _1(x))\;\;=\;\;\rho (\omega (x))\;\;=\;\;\omega (x)+(g'),
\end{eqnarray*}
therefore $\tilde{\omega}$ extends $\omega$ to $\tilde{R}[A]$. If we restrict $\tilde{\omega}$ to $R[A]$, then we have the homomorphism we are looking for, under the conditions that $E'=F'(\alpha ')$ and that $E'$ is a Galois extension of $F'$.\\

\noindent \underline{5. $E' = F'(\alpha ' )$}: As 
$$
\bar{\phi }^{-1}(\alpha )= X + (f)
$$
and
$$
\gamma (X+(f)) = \rho (\hat{\omega}_1(X)) = \rho (X) = X + (g'),
$$
we have
$$
\alpha ' = \tilde{\omega }(\alpha ) = X + (g')
$$
and, by Result \ref{lemSPLIT1},
$$
F'(\alpha ') = F'(X+(g')) = \frac{F'[X]}{(g')} = E'.
$$
\vspace{1mm}

\noindent \underline{6. $E'$ is a Galois extension of $F'$}: As $\charac F'=0$, we only need to show that $E'$ is a normal extension of $F'$. Let $\alpha _1, \ldots ,\alpha _n$ be the roots of $f$. Since $f$ is irreducible over $F$ and $E$ is a splitting field of $f$, Theorem \ref{thGALPOLYirred1} ensures that the Galois group $G=G(E/F)$ acts transitively on the roots of $f$. This implies that the roots of $f$ belong to $A$, because $\alpha \in A$. Moreover, the roots of $f'$ are $\tilde{\omega}(\alpha _1), \ldots , \tilde{\omega}(\alpha _n)$, since, by the relations between the roots of a polynomial and its coefficients, 
$$
f'(X) = (-\tilde{\omega}(\alpha _1)+X)\cdots (-\tilde{\omega}(\alpha _n)+X).
$$
Consequently, 
$$
E'=F'(\alpha ')=F'(\tilde{\omega}(\alpha _1), \ldots ,\tilde{\omega}(\alpha _n))
$$ 
is a splitting field of $f'$ and, by Theorem \ref{NORMALth1}, $E'$ is a normal extension of $F'$. \\

\noindent \underline{7. The special case $g'=f'$}: In this case, $f'$ is irreducible. As above, let $\alpha _1, \ldots ,\alpha _n$ be the conjugates of $\alpha$. Since $E=F(\alpha )$, from Result \ref{propSPLIT2} there exists a unique $\sigma _i\in G$ such that $\sigma _i(\alpha )=\alpha _i$. Similarly, $\alpha _1', \ldots ,\alpha _n'$ are the conjugates of $\alpha '$ and, since $E'=F'(\alpha ')$ and $f'$ is irreducible over $F'$, Result \ref{propSPLIT2} ensures the existence of a unique $\sigma _i'\in G'=Gal(E'/F')$ such that $\sigma _i'(\alpha ')=\alpha _i'$.

From Step $1.$ of our proof (the definition of $u$), the values of $\alpha _1,\ldots ,\alpha _n$ are distinct. Consequently, the automorphisms $\sigma _1,\ldots ,\sigma _n$ are different elements of $G$. Moreover, Result \ref{thGALGRP1} ensures that $G$ has cardinal $n$. Thus $G=\{\sigma _1, \ldots ,\sigma _n\}$. Similarly, since $f'$ is irreducible, $G'$ is of cardinal $n$ and $G'=\{\sigma _1', \ldots ,\sigma _n'\}$.

We now define a mapping $\Phi$ from $G$ into $G'$ by setting $\Phi (\sigma _i)=\sigma _i'$. We will prove that this mapping is an isomorphism. First we will show that
\begin{equation}\label{eqnHILBERTprop1}
\forall s\in \tilde{R}[A], \;\forall \sigma _i\in G, \;\;\tilde{\omega}(\sigma _i(s)) = \sigma _i'(\tilde{\omega}(s)).
\end{equation}
As $\tilde{R}[A]=\tilde{R}[\alpha ]$, it is sufficient to prove the identity for $\alpha$ and for the elements of $\tilde{R}$. For $\alpha$ we have
$$
\tilde{\omega}(\sigma _i(\alpha ) = \tilde{\omega }(\alpha _i) = \alpha _i' = \sigma _i'(\alpha ') = \sigma _i'(\tilde{\omega }(\alpha )).
$$
If $x\in \tilde{R}$, then $x\in F$, hence
\begin{eqnarray*}
\tilde{\omega} (\sigma _i(x)) &=& \tilde{\omega}(x)\;\;= \;\;\gamma \circ\bar{\phi}^{-1}(x) \;\;=\;\;\gamma (x+(f)) \;\;=\;\;\rho\circ \hat{\omega}_1(x)\\
						& =& \hat{\omega}_1(x) +(g')\;\; = \;\;\hat{\omega}_1(x) +(f'),
\end{eqnarray*}
because $g'=f'$. However, $\hat{\omega}_1(x) \in F'$, therefore
$$
\hat{\omega}_1(x) + (f') = \sigma _i'(\hat{\omega}_1(x)+(f')).
$$
Thus
$$
\tilde{\omega}(\sigma _i(x)) = \sigma _i'(\tilde{\omega }(x)).
$$
It follows that the identity (\ref{eqnHILBERTprop1}) applies. We now use this identity to prove that $\Phi$ is a homomorphism. Since $\sigma _i(\alpha )\in A$, for $i=1,\ldots ,n$, we have
\begin{eqnarray*}
(\sigma _i\sigma _j)'(\alpha ') &=& (\sigma _i\sigma _j)'(\tilde{\omega}(\alpha ))\\
						&=& \tilde{\omega}((\sigma _i\sigma _j)(\alpha ))\\
						&=& \tilde{\omega}(\sigma _i(\sigma _j(\alpha )))\\
						&=& \sigma _i'(\tilde{\omega}(\sigma _j(\alpha )))\\
						&=& \sigma _i'\sigma _j'(\tilde{\omega}(\alpha ))\\
						&=& \sigma _i'\sigma _j'(\alpha ').
\end{eqnarray*}
Therefore $\Phi$ is a homomorphism. Clearly $\Phi$ is surjective. As $\vert G\vert = \vert G'\vert$, $\Phi$ is also injective, hence an isomorphism. This finishes the proof.\bx\\

In inverse Galois theory we are confronted with the problem of determining whether a given group $H$ may be considered as the Galois group of a Galois extension $E$ of a certain field $F$, usually $\Q$. This may be difficult to decide directly. However, it may be possible to take another field $F'$ and find a Galois extension $E'$ of this field such that $H$ is isomorphic to the Galois group $Gal(E'/F')$, which is in turn isomorphic to the group $Gal(E/F)$. We will now consider this question.\\

Let $F$ be a Hilbertian field and $E$ a Galois extension of degree $n$ of $F(X)$, the field of fractions of the polynomial ring $F[X]$. From Proposition \ref{proposHILBERTprop1}, we know that there is an element $\alpha \in E$, such that $E=F(X)(\alpha )$ and $f(Y)=m(\alpha ,F(X))\in F[X][Y]$. Then $f$ is irreducible in $F[X][Y]$ (in the sense of our definition at the beginning of the article).

\begin{theorem}\label{thmHILBERTprop1} For an infinite number of values $b\in F$, $f_b(Y)=f(b,Y)$ is irreducible in $F[Y]$ and $E' = F[Y]/(f_b)$ is a Galois extension of $F$, with $G=Gal(E/F(X))$ isomorphic to $G'=Gal(E'/F)$.
\end{theorem}

\noindent \textsc{proof} We apply Lemma \ref{lemHILBERTprop2}. If $A$ is the set of roots of $f$ in $E$, then
$$
\forall x\in A\; \forall \sigma \in G=Gal(E/F(X)), \;\sigma (x) \in A
$$
and $\alpha \in A$. We choose $b\in F$ and consider the homomorphism
$$
\omega _b : F[X] \longrightarrow F, g\longmapsto g(b).
$$
From Lemma \ref{lemHILBERTprop2}, there exists $u\in F[X]$ such that, if $\omega _b(u)\neq 0$, then there is an isomorphism $\Phi _b$ from 
$Gal (E/F(X))$ onto $Gal(E'/F)$, if $f_b$ is irreducible. Indeed, if $f(X,Y) = \sum _{i=0}^na_i(X)Y^i$, then
$$
f'(Y) = \sum _{i=0}^n\omega _b(a_i(X))Y^i = \sum _{i=0}^na_i(b)Y^i = f_b(Y),
$$
where $f'$ is defined as in Lemma \ref{lemHILBERTprop2}. Moreover, in Step $4.$ of Lemma \ref{lemHILBERTprop2}  we saw that $E'=F'[X]/(f')$, which leads to the form of $E'$ in the statement of the theorem. We notice that $u\in F[X]$, so that $u(b)=0$ for a finite number of values $b$. Eliminating these values from the infinite number of values $b$ with $f_b$ irreducible leaves us with an infinite number of values $b$, hence the result.\bx\\

\noindent {\bf Remark} Suppose that we have a finite group $H$ and we wish to know whether it can be represented as a Galois group over a given Hilbertian field $F$, then we may look for a Galois extension $E$ of the field $F(X)$ such that $H$ is isomorphic to the Galois group $G=Gal(E/F(X))$. (As $E$ is a finite extension of a field of functions, we say that $E$ is a function field.) From the theorem, we know that there is a Galois extension $E'$ of $F$ such that $H$ is isomorphic to $G'=Gal(E'/F)$. However, we do not know how to find such an extension.\\

\noindent \underline{\bf 2. A characterization of Hilbertian fields} \\

In this section we continue our discussion of extensions of Hilbertian fields and find a useful characterization of Hilbertian fields. As in Theorem  \ref{thmHILBERTprop1}, we consider a Hilbertian field $F$ and $E$ a Galois extension of degree $n$ of $F(X)$ and we take an element $\alpha \in E$, such that $E=F(X)(\alpha )$ and let $f(X,Y)=m(\alpha ,F(X))\in F[X][Y]$. 

\begin{proposition}\label{propHILBERText1} Let $L$ be a finite extension of $F$ such that $L(X)\subset E$ and $h\in L[X][Y]$ irreducible with roots in $E$. Up to a finite number of the values $b\in F$ such that $f_b(Y)$ is irreducible $h_b(Y)$ is irreducible in $L[Y]$.
\end{proposition}

\noindent \textsc{proof} Let $\beta _1', \ldots ,\beta _m'$ be the roots of $h_b$ in an extension $E'$ of $L$, where $f_b$ is irreducible in $F[Y]\subset L[Y]$. If $h_b$ is reducible, then we can write $h_b=u_bv_b$, with $u_b,v_b\in L[Y]$ and $\deg u_b>0$, $\deg v_b>0$. Without loss of generality we may write 
\begin{eqnarray*}
u_b(Y) &=& \gamma (-\beta _1'+Y)\cdots (-\beta _s'+Y)\\
v_b(Y) &=& \delta (-\beta _{s+1}'+Y)\cdots (-\beta _m'+Y),
\end{eqnarray*}
with $\gamma ,\delta \in L$. If $\sigma \in Gal (E'/L)$, then $\sigma$ permutes the roots of $u_b$ and of $v_b$, which cannot be possible if $Gal (E'/L)$ acts transitively on the roots. It follows that, if $Gal (E'/L)$ acts transitively on the roots of $h_b$, then $h_b$ is irreducible. We will show that, with the exception of a finite set of values of $b$, there is an extension $E'$ of $L$ containing the roots of $h_b$ and such that the Galois group $Gal(E'/L)$ acts transitively on these roots. 

As $L$ is a finite extension of $F$ and $\charac F=0$, by the primitive element theorem, there exists $x\in L$ such that $L=F(x)$. We notice that
$$
F\subset F(x) = L \Longrightarrow F(X)\subset L(X) \subset E.
$$
In addition, as $x\in L$, we have $x\in L(X)$, so we may consider the conjugates of $x$ over $F(X)$. The roots $\beta _1, \ldots ,\beta _m$ of $h(X,Y)$ are, by hypothesis, in $E$. We let $A$ be the subset of $E$ composed of $x$, with its conjugates over $F(X)$, $\alpha$, with its conjugates over $F(X)$, and the roots of $h(X,Y)$. If $\sigma \in Gal(E/F(X))$ and $a\in A$, then clearly $\sigma (a)\in A$. We notice that the roots of $h_b$ are distinct if the discriminant $\Delta (h_b)\neq 0$. However, $\Delta (h_b)=\Delta (h)(b)$, which has the value $0$ for a finite number of values of $b$. Therefore, if we exclude these values, we can be sure that the roots of $h_b$ are distinct.

As in Theorem \ref{thmHILBERTprop1}, we consider the valuation homomorphism 
$$
\omega _b : F[X] \longrightarrow F, g\longmapsto g(b),
$$
where $b\in F$. Lemma \ref{lemHILBERTprop2} ensures the existence of $u\in F[X]$ such that, if $\omega _b(u)\neq 0$, we may extend $\omega _b$ to a homomorphism 
$$
\tilde{\omega}_b : F[X][A] \longrightarrow E',
$$
where $E'$ is a Galois extension of $F$. Moreover, $u(b)=0$ for a finite number of values $b$, so we may suppose that $\omega _b(u)\neq 0$. As a generator of $L$ over $F$ is included in $A$, we have $L\subset F[X][A]$. Also, $\tilde{\omega}_b$ restricted to $F$ is the identity and so $\tilde{\omega}_b$ restricted to $L$ is not trivial. Given that a ring homomorphism of a field into another field is either trivial or injective, it must be so that $\tilde{\omega}_b$ restricted to $L$ is injective. We will note $\tilde{L}$ the image of $L$ under $\tilde{\omega}_b$ in $E'$. As $\tilde{L}$ is isomorphic to $L$, $E'$ is an extension of $L$. 

We now let $\chi _b$ be the natural homomorphism from $F[X][A][Y]$ into $E'[Y]$ generated by $\tilde{\omega}_b$. We may identify $F[A][X]$ and $F[X][A]$. To see this it is sufficient to notice that, If $A=\{a_1, \ldots ,a_n\}$, then $F[A][X]$ is composed of sums of expressions of the form $ca_1^{s_1}\ldots a_n^{s_n}X^t$ and $F[X][A]$ of sums of expressions of the form $cX^ta_1^{s_1}\ldots a_n^{s_n}$, where $c\in F$ and $s_1, \ldots s_n, t\in \N$. Thus we may consider that $L[X]=F(x)[X]$ is included in $F[X][A]$, because $x\in A$. We may write
$$
h(X,Y) = \sum _{i=0}^mh_i(X)Y^i,
$$
where the coefficients $h_i(X)$ lie in $L[X]$, a subset of $F[X][A]$. We are interested in finding the form of $\chi _b(h)$. We may write 
$$
h_i(X) = \sum _{j=1}^na_{ij}X^j,
$$
with $a_{ij}\in L$. As $L=F(x)$, we can express each $a_{ij}$ in the form
$$
a_{ij} = \sum _{k=1}^mu_k^{ij}x^k,
$$
with $u_k^{ij}\in F$. Hence,
\begin{eqnarray*}
\tilde{\omega}_b\left(\sum _{j=1}^na_{ij}X^j\right) &=& \sum _{j=1}^n\tilde{\omega}_b(a_{ij})\tilde{\omega}_b(X^j)\\
								&=& \sum _{j=1}^n\tilde{\omega}_b\left(\sum _{k=1}^mu_k^{ij}x^k\right)b^j\\
								&=& \sum _{j=1}^n\left(\sum _{k=1}^mu_k^{ij}\tilde{\omega}_b(x)^k\right)b^j.
\end{eqnarray*}
Therefore, if we identify $\tilde{L}$ and $L$ and consider the polynomial $\tilde{h}\in \tilde{L}[X][Y]$ corresponding to $h\in L[X][Y]$, then we find an expression for the coefficients of $\tilde{h}_b$ by replacing $x$ with $\tilde{\omega}_b(x)$. Let us write $\tilde{h}_b$ for this polynomial. Then we have
$$
\chi _b (h(X,Y) )= \tilde{h}_b(Y).
$$
We notice that the roots of $h$ lie in $A$ so their images under $\tilde{\omega}_b$ are in $E'$. We have
$$
h(X,Y) = h_m(X)(-\beta _1+Y)\cdots (-\beta _m+Y)\Longrightarrow \tilde{h}_b(Y) = \tilde{h}_m(b)(-\beta _1'+Y)\cdots (-\beta _m'+Y),
$$
where $\beta j' = \chi _b(\beta _i)$, for some $i$.  So we have found expressions for the roots of $\tilde{h}_b$.

From our work at the beginning of the proof, to show that $\tilde{h}_b$ is irreducible, it is sufficient to prove that the Galois group $Gal(E'/\tilde{L})$ acts transitively on the roots $\beta _j'$ (for any $b$ in the infinite set we have retained). This we will now do. As $h(X,Y)$ is irreducible in $L[X][Y]$, from Gauss's lemma, $h(X,Y)$ is also irreducible in $L(X)[Y]$. $E$ is a normal extension of $F(X)$ containing $L(X)$, hence, by Result \ref{NORMALprop2}, $E$ is a normal extension of $L(X)$. Moreover, by Result \ref{NORMALth1}, $E$ is the splitting field of a polynomial $g\in L(X)[Y]$. From Result \ref{thGALPOLYirred1}, the Galois group $G_1=Gal(\bar{E}/L(X))$, where $\bar{E}$ is a splitting field of $h$ included in $E$, acts transitively on the roots $\beta _i$. Now, Result \ref{thSPLIT2}, with $F=F'=\bar{E}$ and $f=f^*=g$, implies that any element $\sigma \in G_1$ may be extended to an element $\tilde{\sigma}\in G_2=Gal (E/L(X))$. This implies that $G_2$ acts transitively on the roots $\beta _i$. Supposing that $f_b$ is irreducible, then, from Theorem \ref{thmHILBERTprop1}, $G=Gal(E/F(X))$ is isomorphic to $G'=Gal(E'/F)$, where $E'=F[Y]/(f_b)$. (This is the same $E'$ as that obtained earlier in the proof after applying Lemma \ref{lemHILBERTprop2}.) 

If we restrict the isomorphism from $G$ onto $G'$ to the subgroup $G_2$, then we obtain a subgroup $G_2'$ of $G'$. We claim that $G_2'$ is a subgroup of the Galois group $G''=Gal(E'/\tilde{L})$. To prove this we need to show that the automorphisms of $G_2'$ fix the elements of $\tilde{L}$. We use the identity (\ref{eqnHILBERTprop1}) and the explicit form of the isomorphism $\Phi$ from $G$ onto $G'$. If $\tilde{z}\in \tilde{L}$, $\sigma \in G_2$ and $\sigma '$ the corresponding automorphism in $G_2'$, then
$$
\sigma '(\tilde{z}) = \sigma '(\tilde{\omega }_b(z)) = \tilde{\omega}_b(\sigma (z)) =\tilde{\omega}_b(z) = \tilde{z}.
$$
This proves that $G_2'\subset G''$ and so the claim.

We now show that  $G''$ acts transitively on the roots $\beta _j'$. There exists $\sigma \in G_2$ such that $\sigma (\beta _i)=\beta _j$, because $G_2$ acts transitively on the roots $\beta _i$. Let $\sigma '$ be the element of $G_2'$ corresponding to $\sigma$. Then, using the identity (\ref{eqnHILBERTprop1}) again, we have
$$
\sigma '(\beta _i') = \sigma'(\tilde{\omega }_b(\beta _i)) = \tilde{\omega}_b(\sigma (\beta _i)) = \tilde{\omega }_b(\beta _j) = \beta _j'.
$$
Therefore $G''$ acts transitively on the roots $\beta _j'$. This finishes the proof.\bx\\

We can now establish the characterization of Hilbert fields, which we referred to at the beginning of the section.

\begin{theorem}\label{thHILBERText1} The field $F$ is Hilbertian if and only if, for every finite extension $L$ of $F$ and finite set of irreducible polynomials $h_1, \ldots ,h_k\in L[X,Y]$, there is an infinite number of values $b\in F$ such that $h_{i,b}(Y)=h_i(b,Y)\in L[Y]$ is irreducible, for $i=1,\ldots ,k$.
\end{theorem}

\noindent \textsc{proof} Let $F$ be an Hilbertian field, $L$ a finite extension of $F$ and $h_1, \ldots ,h_k$ irreducible polynomials in $L[X,Y]=L[X][Y]$. (From Gauss's Lemma, these polynomials are irreducible in $L(X)[Y]$.) Adding the roots of the polynomials $h_i$ to $L(X)$, we obtain a finite extension $M$ of $L(X)$. From Result \ref{thSEPext4}, $L(X)$ is a finite extension of $F(X)$, so $M$ is a finite extension of $F(X)$. Now let $E$ be a normal closure of $M$ over $F(X)$. Then $E$ is finite Galois extension of $F(X)$ (Result \ref{NormalClos1}). From Result \ref{proposHILBERTprop1}, we may find $\alpha \in E$ such that $E=F(\alpha )$ and $m(\alpha ,F(X)) \in F[X][Y]$. As usual, we write $f(X)$ for $m(\alpha ,F(X))$. In addition, $L(X)\subset E$ and the roots of each $h_i$ are in $E$. From Proposition \ref{propHILBERText1}, for each $i$, for all but a finite number of the $b\in F$ such that $f_b(Y)$ is irreducible, $h_{i,b}(Y)$ is irreducible. It follows that $h_{i,b}(Y)$ is irreducible, for all $i$, for an infinite number of values of $b$. 

The converse is elementary. We only need to choose $L=F$ and $k=1$.\bx\\

We know that $\Q$ is a Hilbertian field. The next two results will show us that there are many other such fields.

\begin{theorem} If $F$ is a Hilbertian field, then any finite extension $E$ of $F$ is also Hilbertian.
\end{theorem}

\noindent \textsc{proof} Let $E$ be a finite extension of the Hilbertian field $F$ and $f\in E[X,Y]$ irreducible. From Theorem \ref{thHILBERText1}, there are infinitely many $b\in F$, such that $f_b$ is irreducible. As $F\subset E$, there are infinitely many $b\in E$ such that $f_b$ is irreducible. Hence $E$ is Hilbertian.\bx\\

We recall that a number field is a finite extension in $\C$ of the field $\Q$. 

\begin{corollary} Number fields are Hilbertian.
\end{corollary}

\vspace{0.2cm}

\noindent \underline{\bf 3. The Kronecker specialization} \\

We have considered extensions of fields of fractions of polynomials in one variable. In this section we aim to consider the case of polynomials in several variables. We use a tool known as the Kronecker specialization. We fix an integer $d>1$. For a field $F$ and an integer $k>2$, we define a mapping $S_d$ from $F[X_1,\ldots ,X_k]$ into $F[X,Y]$ by
$$
S_d(f)(X,Y) = f(X,Y,Y^d, Y^{d^2}, \ldots ,Y^{d^{k-2}}).
$$ 
The mapping $S_d$ is called a Kronecker specialization. We note $V_d$ the collection of polynomials in $F[X_1, \ldots ,X_k]$ whose degree is less than $d$ in each variable $X_2, \ldots ,X_k]$ and $W_d$ the collection of polynomials in $F[X,Y]$ whose degree is inferior to $d^{k-1}$ in the variable $Y$. 

\begin{proposition} The mapping $S_d$ defines a bijection from $V_d$ onto $W_d$. 
\end{proposition}

\noindent \textsc{proof} First let $f\in V_d$ be a monomial. Then $f$ has the form $aX_1^{\alpha _1}\cdots X_k^{\alpha _k}$ and so
$$
S_d(f) = aX^{\alpha _1}Y^{\alpha _2+\alpha _3d + \alpha _4d^2 +\ldots +\alpha _kd^{k-2}}.
$$
If $g\in V_d$, with $g=bX_1^{\beta _1}\cdots X_k^{\beta _k}$, and $S_d(g)=S_d(f)$, then 
$$
b=a \qquad \beta _1= \alpha _1\qquad \beta _2+\beta _3d + \beta _4d^2 +\ldots +\beta _kd^{k-2}=\alpha _2+\alpha _3d + \alpha _4d^2 +\ldots +\alpha _kd^{k-2}.
$$
As the representation in base $d$ is unique, we must have $\alpha _i=\beta _i$, for all $i\geq 2$. Hence $g=f$ and $S_d$ defines an injection from the monomials in $V_d$ into the monomials in $W_d$. Also, any integer $s<d^{k-1}$ has a unique representation in base $d$:
$$
s = \alpha _2 + \alpha _3d + \alpha _4d^2 + \cdots + \alpha _kd^{k-2},
$$
which implies that $S_d$ restricted to the monomials of $V_d$ defines a bijection onto the monomials of $W_d$. As 
$$
S_d(m_1 + \cdots + m_k) = S_d(m_1) + \cdots + S_d(m_k),
$$
for monomials $m_1, \ldots ,m_k$, the mapping $S_d$ defines a bijection from $V_d$ onto $W_d$\bx\\

\noindent {\bf Remark} It  is not difficult to see that, if the product $fg$ is in $V_d$, then 
$$
S_d(fg) = S_d(f)S_d(g).
$$

We now see that the Hilbertian property "carries over" to multivariable polynomials. (We advise the reader to return to the beginning of the article to revise the definition of an irreducible polynomial in several variables.)

\begin{theorem}\label{thHILBERTkron1} If $F$ is a Hilbertian field and $f\in F[X_1, \ldots ,X_k]$ is irreducible, then there exists an infinite number of values $b\in F$ such that $f(b,X_2, \ldots ,X_k)\in F[X_2, \ldots ,X_k]$ is irreducible.
\end{theorem}

\noindent \textsc{proof} Let $d$ be an integer superior to the degree of each variable $X_i$ in $f$. We can write
\begin{equation}\label{eqnHILBERTkron1}
S_d(f)(X,Y) = g(X)\prod _{i\in C}g_i(X,Y),
\end{equation}
where the $g_i$ are irreducible polynomials in $F[X][Y]$ and $C$ is a finite index set. As $F$ is Hilbertian, from Theorem \ref{thHILBERText1}, for an infinite number of values of $b\in F$, $g_{i,b}(Y)=g_i(b,Y)$ is irreducible in $F[Y]$. If we exclude the values $b$ which are roots of $g$, then we still have an infinite number of values of $b\in F$ such that $g_i(b,Y)$ is irreducible, for all $i\in C$. For any $b$ in the remaining set, we have the factorization into prime factors in $F[Y]$:
$$
S_d(f)(b,Y) = g(b)\prod _{i\in C}g_i(b,Y).
$$

For any such $b$, suppose now that $f_b=f(b,X_2,\ldots ,X_k)$ is reducible, i.e., $f_b=uv$, with $u$ and $v$ nonconstant. We may consider $f_b$, $u$ and $v$ as members of $F[X_1,\ldots ,X_k]$. As $f_b\in V_d$, we have $u,v\in V_d$ and
$$
S_d(u)S_d(v) = S_d(uv) = S_d(f_b) = S_d(f)(b,Y) = g(b)\prod _{i\in C}g_i(b,Y).
$$
We now have a partition $\{A,B\}$ of $C$, with $A$ and $B$ nonempty, since $u$ and $v$ are nonconstant, and $\alpha , \beta\in F$ such that $g(b)=\alpha\beta$ and
$$
S_d(u) = \alpha\prod _{i\in A}g_i(b,Y) \qquad \textrm{and}\qquad S_d(v) = \beta\prod _{i\in B}g_i(b,Y).
$$
At this point we set
$$
U(X,Y) = \prod _{i\in A}g_i(X,Y) \qquad \textrm{and}\qquad V(X,Y) = \prod _{i\in B}g_i(X,Y).
$$
As $U,V\in W_d$, there exist unique $\tilde{u}, \tilde{v}\in V_d$ such that $S_d(\tilde{u})=U$ and $S_d(\tilde{v})=V$. Then
$$
S_d(\tilde{u}_b) = S_d(\tilde{u})(b,Y) = U(b,Y) = \prod _{i\in A}g_i(b,Y) = \alpha ^{-1}S_d(u) = S_d(\alpha ^{-1}u),
$$
which implies that $\tilde{u}_b=\alpha ^{-1}u$. In the same way, $\tilde{v}_b=\beta ^{-1}v$ and so
$$
\tilde{u}_b\tilde{v}_b=\alpha ^{-1}\beta ^{-1}uv = g(b)^{-1}f_b.
$$

Our next step is to show that $\tilde{u}\tilde{v}\notin V_d$. If this is not the case, then
$$
S_d(g\tilde{u}\tilde{v}) = gS_d(\tilde{u})S_d(\tilde{v}) = gUV = S_d(f),
$$
from which we deduce that
$$
g\tilde{u}\tilde{v} = f,
$$
which contradicts the irreducibility of $f$, since $\tilde{u}$ and $\tilde{v}$ are nonconstant. Hence $\tilde{u}\tilde{v}\notin V_d$.

We are now in a position to prove that $f_b$ is irreducible for an infinite number of values of $b$. We may consider $f$ as a polynomial in $F[X_1][X_2, \ldots ,X_k]$. Thus a monomial is of the form $a(X_1)X^{n_2}\cdots X^{n_k}$. As $\tilde{u}\tilde{v}\notin V_d$ and $\tilde{u}_b\tilde{v}_b=g(b)^{-1}f_b$, to avoid a contradiction, $b$ must be a root of $a(X_1)$, whenever there is an $i$ such that the power $n_i$ of $X_i$ is greater than $d-1$. We can eliminate these values of $b$ for each such monomial. As the number of these monomials is finite, we are left with an infinite number of values of $b$ for which $f_b$ is irreducible.\bx

\begin{corollary} Let $F$ be a Hilbertian field and $f\in F[X_1, \ldots ,X_k]$, with $k\geq 2$, irreducible. For every polynomial $p\in F[X_1, \ldots ,X_{k-1}]$, there exist elements $b_1, \ldots ,b_{k-1}\in F$ such that $p(b_1, \ldots ,b_{k-1})\neq 0$ and the polynomial $f(b_1, \ldots ,b_{k-1},X_k)$ is irreducible.
\end{corollary}

\noindent \textsc{proof} Let $f\in F[X_1, \ldots ,X_k]$ be irreducible. We aim to prove by induction on $n$ that, for every polynomial $p\in F[X_1, \ldots ,X_n]$, with $n<k$, there exist $b_1,\ldots , b_n\in F$ such that $p(b_1, \ldots ,b_n)\neq 0$ and $f(b_1, \ldots ,b_n, X_{n+1}, \ldots ,X_k)$ is irreducible. First, let $p\in F[X_1]$. From Theorem \ref{thHILBERTkron1}, there is an infinite number of values of $b\in F$ such that $f(b, X_2, \ldots ,X_k)$ is irreducible. If we take one such $b$ which is not a root of $p$, then we have a value of $b$ satisfying the required conditions. Thus the result is true for $n=1$. 

We now suppose that the result is true for $n<k-1$ and consider the case $n+1$. Let $p\in F[X_1, \ldots ,X_{n+1}]$. We may consider $p$ as an element of $F[X_{n+1}][X_1,\ldots ,X_n]$. Each coefficient has a finite number of roots. If $c$ is not one of these roots, then $p(X_1, \ldots ,X_n, c)$ is a nonzero polynomial in $F[X_1, \ldots ,X_n]$. From the induction hypothesis, there exist $b_1, \ldots ,b_n$ such that $p(b_1, \ldots b_n,c)\neq 0$ and $f(b_1, \ldots , b_n, X_{n+1}, \ldots , X_k)$ is irreducible. Using Theorem \ref{thHILBERTkron1}, we know that there is an infinite number of values of $b$, such that $f(b_1, \ldots ,b_n, b,X_{n+2}, \ldots ,X_k)$ is irreducible. As $p(b_1, \ldots ,b_n,X_{n+1})\in F[X_{n+1}]$, if we eliminate any $b$ which is a root of this polynomial, then we have $p(b_1, \ldots ,b_n, b)\neq 0$. Therefore the result is true for $n+1$. This finishes the induction step. \bx\\

%We recall that, if $R$ is a unique factorization domain and $f\in F[X]$, then the {\it content}\index{content of the polynomial} $f$, which we write $c(f)$, is the pgcd of the coefficients of $f$. We say that a polynomial is {\it primitive}\index{primitive polynomial} if its content is $1$. Clearly, we may write $f=c(f)h$, where the $c(h)=1$. 

To prove the next theorem, which will provide us with more Hilbertian fields, we will need Result \ref{ufd}. 

%Let $R$ be a unique factorization domain, with quotient field $F$, and $f\in R[X]$. Then, if $f$ is nonconstant and irreducible in $R[X]$, then $f$ is irreducible in $F[X]$. On the other hand, if $f$ is primitive and irreducible in $F[X]$, then $f$ is irreducible in $R[X]$. 

\begin{theorem}\label{thHILBERTkron2} If $F$ is Hilbertian field, then $F(X_1, \ldots , X_k)$ is Hilbertian field, for any $k\in \N^*$.
\end{theorem}

\noindent \textsc{proof} Let $f\in F(X_1, \ldots ,X_k)[X,Y]$ be irreducible. There exists $g\in F(X_1, \ldots ,X_k)^*$ such that $gf \in F[X_1, \ldots X_k,X][Y]$. We may write $gf=c(gf)h$, where $c(h)=1$. We have
$$
h = \frac{g}{c(gf)}f \qquad \textrm{and}\qquad \frac{g}{c(gf)}\in F(X_1, \ldots ,X_k,X)^*.
$$
We notice that $h\in F(X_1,\ldots ,X_k,X)[Y]$ is primitive and irreducible. From Gauss's Lemma, $h$ is also irreducible in $F[X_1,\ldots ,X_k,X][Y]$. Now, using Theorem \ref{thHILBERTkron1}, we see that there are infinitely many values of $b\in R$ such that $h(X_1,\ldots ,X_k,b,Y)$ is irreducible in $F[X_1,\ldots ,X_k][Y]$. Using Gauss's Lemma again, we see that $h(X_1,\ldots ,X_k,b,Y)$ is irreducible in $F(X_1,\ldots ,X_k)[Y]$. To finish, we notice that
$$
f(X_1, \ldots ,X_k,b,Y) = \frac{c(gf)_b}{g}h(X_1,\ldots ,X_k,b,Y),
$$
where $c(gf)_b$ is the polynomial in $F[X_1,\ldots ,X_k]$ obtained by replacing the variable $X$ by $b$. Now $c(fg)\in F(X_1,\ldots X_k)[X]$, so there can only be a finite number of values of $b$ such that $c(gf)_b=0$. If we exclude these values, then $\frac{c(gf)_b}{g}$ is a unit in $F(X_1,\ldots ,X_k)$ and so $f(X_1,\ldots ,X_k,b,Y)$ is irreducible in $F(X_1,\ldots X_k)[Y]$. We have shown that $F(X_1, \ldots ,X_k)$ is Hilbertian.\bx\\

We may now extend Theorem \ref{thmHILBERTprop1}.

\begin{theorem}\label{HREDthm1a} If $F$ is an Hilbertian field and $E$ a Galois extension of $F(X_1, \ldots ,X_k)$, then there exists a Galois extension $E'$ of $F$ such that $Gal (E/F(X_1, \ldots ,X_k))$ is isomorphic to $Gal(E'/F)$.
\end{theorem}

\noindent \textsc{proof} We prove this result by induction on $k$. For $k=1$, it is sufficient to apply Theorem \ref{thmHILBERTprop1}. Suppose now that the result is true for $k$ and let us consider the case $k+1$. $E$ is a Galois extension of $F(X_1, \ldots ,X_k,X_{k+1})$. From  Theorem \ref{thHILBERTkron2}, $F(X_1, \ldots ,X_k)$ is Hilbertian and so there is a Galois extension $E'$ of $F(X_1, \ldots ,X_k)$ such that 
$$
Gal(E/F(X_1, \ldots ,X_k, X_{k+1})) \simeq Gal (E'/F(X_1, \ldots ,X_k)).
$$
From the induction hypothesis there is an extension $E''$ such that 
$$
Gal (E'/F(X_1, \ldots ,X_k))\simeq Gal (E''/F),
$$
therefore 
$$
Gal(E/F(X_1, \ldots ,X_k, X_{k+1}))\simeq Gal (E''/F).
$$
This finishes the induction step and the proof.\bx\\

\noindent {\bf Remark} We recall that the general polynomial of degree $k$ over a field $F$ has the form
$$
f(Y) = Y^n -X_1Y^{n-1} + X_2Y^{n-2} + \cdots +(-1)^{n-1}X_{n-1}Y + (-1)^nX_n \in F(X_1, \ldots ,X_k)[Y],
$$
where $F(X_1, \ldots ,X_k)$ is the rational function field over $F$ in $k$ variables. The Galois group of $f$ is the symmetric group $S_k$ (see \cite{spindler}). In particular, this is the case if $F=\Q$. Writing $E$ for the splitting field of $f$ over $\Q(X_1, \ldots ,X_k)$ and using the fact that $\Q$ is a Hilbertian field, from Theorem \ref{HREDthm1a} we deduce that $S_k$ is realizable as a Galois group over $\Q$. \\

\vspace{0.2cm}

\noindent \underline{4. \bf Proof of Hilbert's irreducibility theorem} \\

Above we stated without proof Hilbert's irreducibility theorem, namely that the field of rational numbers $\Q$ is a Hilbertian field, i.e., for any $f\in \Q[X,Y]$, which is irreducible, there exists an infinite number of values $b\in\Q$ such that $f(b,Y)\in \Q[Y]$ is irreducible. We aim now to provide a detailed proof of this result. Our proof is based on that given in \cite{hadlock}, with modifications. We now suppose that our field is $\C$, the complex numbers. For any $z\in \C$, the polynomial $f(z, Y)$ has $n$ roots which may or may not be distinct. For any $z$, we may write $u_1(z), \ldots ,u_n(z)$ for these roots. Thus we obtain $n$ functions defined on $\C$. If $b$ is a regular value, then the roots $u_1(b), \ldots ,u_n(b)$ are distinct. We can say more.

\begin{lemma} \label{lemHILBIRRED1} If $b$ is a regular value of the polynomial $f(X,Y)$ over $\C$, then there is a neighbourhood $N$ of $b$ in $\C$ such that the functions $u_1(z), \ldots ,u_n(z)$ are analytic on $N$.
\end{lemma} 

\noindent \textsc{proof} We will suppose that the functions $u_i$ exist and find the possible forms, then we will show that the functions so obtained satisfy the conditions. Let $u(z)$ be one such function and, to begin, we will suppose that $b=0$ and $u(0)=0$. We seek a power series
$$
u(z) = \sum _{k=1}^{\infty}b_kz^k,
$$
with $b_k\in \C$, which converges on some set $N=\{z\in \C: \vert z\vert <R\}$, with $R>0$, within which
$$
f(z,u(z)) = 0.
$$
We may write $f(z,u(z))$ in the form
$$
f(z,u) = a_{10}z + a_{01}u + \sum _{i+j\geq 2}a_{ij}z^iu^j.
$$
The sum on the right is finite, because $f$ is a polynomial. Let $f_u$ be the derivative of $f$ with respect to the variable $u$. As $z=0$ is a regular value, $f_u(0,0) = a_{01}\neq 0$, so we may write
\begin{equation}\label{eqnHILBERTreduc1}
f(z,u) = -a_{01}\left(-\frac{a_{10}}{a_{01}}z -u + \sum _{i+j\geq 2}-\frac{a_{ij}}{a_{01}}z^iu^j\right) = -a_{01}\left(a_{10}'z -u + \sum _{i+j\geq 2}a_{ij}'z^iu^j\right).
\end{equation}
From this we deduce 
$$
f_u(z,u) = -a_{01}\left(-1+g(z,u)\right),
$$
where every monomial of $g$ has degree at least $1$.

We now substitute $u(z)=\sum _{k=1}^{\infty}b_kz^k$ in the equation $f(z,u(z))=0$:
$$
f(z,u(z)) = -a_{01}\left(a_{10}'z - \sum _{k\geq 1}b_kz^k + \sum _{i+j\geq 2}a_{ij}'z^i\left(\sum _{k\geq 1}b_kz^k\right)^j\right) = 0.
$$
As this is a power series in $z$, which we suppose convergent on a neighbourhood of $0$, the coefficient of $z$, namely $a_{10}'-b_1$, has the value $0$, which implies that $b_1=a_{10}'$. 

Our next step is to take a Taylor expansion of the polynomial $f(z,u)$ in $u$, around the point $u_0=\sum _{i=1}^{k-1} b_iz^i$. For $u(z)=\sum _{k\geq 1}b_kz^k$, we obtain
\begin{eqnarray*}
f(z,\sum _{i=1}^{\infty}b_iz^i) &=& f(z, \sum _{i=1}^{k-1}b_ iz^i) +\left(-a_{01}\left(-1+g(z,\sum _{i=1}^{k-1}b_ iz^i)\right)\right)\sum _{i= k}^{\infty}b_iz^i\\
&& +\; \textrm{terms of degree at least $2k$}.
\end{eqnarray*}
The value of this expression is $0$, so the coefficient of each power of $z$ must be $0$. In particular, the coefficient of $z^k$ is $a_{01}b_k$ plus the coefficient $c_k$ of $z^k$ in the expression of the polynomial $f(z, \sum _{i=1}^{k-1}b_ iz^i)$. Therefore $b_k=-\frac{c_k}{a_{01}}$. With this rule and the initial value $b_1=a_{10}'$ the entire sequence of $b_i$s is determined. We now have a candidate for the series $\sum _{k=1}^{\infty} b_kz^k$. We must show that it has a positive radius of convergence. We will do this by constructing a power series $\sum _{k=1}^{\infty}A_kz^k$ with positive radius of convergence $R$ and such that, for every $k$, $A_k\geq \vert b_k\vert$. 

For each $k$, there is a polynomial in several variables with positive integer coefficients, which we note $p_k$, such that $b_k$ is is the value of $p_k$ evaluated at the set of coefficients $a_{ij}'$. We will write $b_k=p_k(a_{ij}')$. Let $A$ be an element of $\Z$ such that $A\geq \vert a_{ij}'\vert$, for all $a_{ij}'$. If we replace $a_{ij}'$ by $A$ in the equation $(\ref{eqnHILBERTreduc1})$, then we obtain
$$
h(z,v) = Az - v + A\sum _{i+j\geq 2}z^iv ^j = 0.
$$
Now, redoing the calculations which we did on the equation $(\ref{eqnHILBERTreduc1})$, we obtain a solution of the form $\sum _{k=1}^{\infty}A_kz^k$, where $A_k$ is the polynomial expression derived from $p_k(a_{ij}')$ by replacing the $a_{ij}'$s by $A$. Clearly $A_k\geq \vert b_k\vert$, for every $k$, since the coefficients of the $p_k$ are positive. For $\vert z\vert <1$ and $\vert u\vert <1$, we have
\begin{eqnarray*}
0 &=& Az - v + Az^0\sum _{j=2}^{\infty} v^j + Az ^1\sum _{j=1}^{\infty}v^j + A\sum _{i=2}^{\infty}z^i\left(\sum _{j=0}^{\infty}v^j\right)\\
&=& Az - v + A\frac{v^2}{1-v} + A\frac{zv}{1-v} + A\left(\frac{z^2}{1-t}\frac{1}{1-z}\right).
\end{eqnarray*}
Multiplying by $1-v$ we obtain
\begin{equation}\label{eqnHILBERTreduc2}
0 = Az(1-v) -v(1-v) + Av^2 + Azv + A\frac{z^2}{1-z} = (A+1)v^2 - v +\frac{Az}{1-z},
\end{equation}
which is a quadratic equation in $v$. If we set, for $z$ sufficiently small, 
$$
v(z) = \frac{1-\left(1-4(A+1)\frac{Az}{1-z}\right)^{\frac{1}{2}}}{2(A+1)},
$$
then $v(z)$ is the solution of the equation $(\ref{eqnHILBERTreduc2})$, with $v(0)=0$, where the square root in the expression is the principal value (for example, see \cite{henrici}). We would like to show that $v(z)$ is analytic for $z$ sufficiently small. However, the principal value of the function $w_1(z)=(1-z)^ {\frac{1}{2}}$ is analytic, which is also the case for the function $w_2(z)=\frac{z}{1-z}$. Moreover, for $z$ sufficiently small, $\vert w_2(z)\vert <1$. Therefore, given that the composition of analytic functions is analytic, we see that $v(z)$ is analytic in a neighbourhood of $0$. It now follows that $u(z)$ is analytic in a neighbourhood of $0$.

We have supposed that $u(0)= 0$. If we have $u(0)=a$, where $a$ is not necessarily $0$, then, writing $y(z)=u(z)-a$, we have $y(0)=0$ and so we have $y(z)$ analytic for $z$ sufficiently small, which implies that $u(z)$ is analytic for $z$ sufficiently small. Finally, if we suppose that $u(b)=a$, then writing $y(z)=u(b-z)$, we have $y(0)=a$ and $y(z)$ is analytic for $z$ sufficiently close to $0$, which implies that $u(z)$ is analytic for $u$ sufficiently close to $b$.

For each root $a_i$ of $f_b$ we can find an analytic root function $u_i(z)$, with $u_i(b)=a_i$, defined on a neighbourhood $N_i$ of $b$. As the $a_i$ are distincts, by continuity, we may choose a neighbourhood $N$ of $b$ such that the functions $u_i(z)$ are analytic on $N$ and, for  any $z\in N$, distinct. This finishes the proof.\bx\\

\noindent {\bf Remark} For $z\in N$, we may write
$$
f(z,Y) = a_n(z)\prod _{i=1}^n \left(-u_i(z) + Y\right),
$$
where $a_n(z)$ is a polynomial in $z$ and the $u_n$ are analytic functions defined on $N$.\\

We need another preliminary result. We take $m+1$ increasing values of the real variable $t$: $t_0<t_1<t_2<\cdots <t_m$ and we write $V_m$ for the Vandermonde determinant of the $t_i$s, i.e.,
\begin{equation*}
V_m = \begin{vmatrix}
1 & t_0 & t_0^2 & \ldots & t_0^{m-1} & t_0^m\\
1 & t_1 & t_1^2 & \ldots & t_1^{m-1} & t_1^m\\ 
\vdots & \vdots & \vdots & \ddots & \vdots & \vdots \\
1 & t_m & t_m^2 & \ldots & t_m^{m-1} & t_m^m
\end{vmatrix}.
\end{equation*}
Now let $f:[t_0,t_m]\longrightarrow \R$ be an $m$ times differentiable function. We set
\begin{equation*}
W_m = \begin{vmatrix}
1 & t_0 & t_0^2 & \ldots & t_0^{m-1} & f(t_0)\\
1 & t_1 & t_1^2 & \ldots & t_1^{m-1} & f(t_1)\\ 
\vdots & \vdots & \vdots & \ddots & \vdots & \vdots \\
1 & t_m & t_m^2 & \ldots & t_m^{m-1} & f(t_m)
\end{vmatrix}.
\end{equation*}

\begin{lemma}\label{lemHILBIRRED2} There exist $u\in (t_0,t_m)$ such that 
$$
\frac{W_m}{V_m} = \frac{f^{(m)}(u)}{m!}.
$$
\end{lemma}

\noindent \textsc{proof} Suppose that $g:[t_0,t_m]\longrightarrow \R$ is an $m$ times differentiable function and that $g(t_i)=f(t_i)$, for all $i$. By Rolle's theorem, in each interval $(t_i, t_{i+1})$, there exists $a_i$ such that $g^{(1)}(a_i)=f^{(1)}(a_i)$. So we have $m$ points $a_0<a_1<\cdots <a_{m-1}$ such that $g^{(1)}(a_i)=f^{(1)}(a_i)$, for all $i$. We now applly Rolle's theorem again and obtain $m-1$ points $b_0<b_1<\cdots <b_{m-2}$ such that $g^{(2)}(b_i)=f^{(2)}(b_i)$, for all $i$. Continuing in the same way, we finally obtain a point $u$ such that $g^{(m)}(u)=f^{(m)}(u)$

We now consider the system
\begin{equation*}
\begin{pmatrix}
1 & t_0 & \ldots & t_0^m\\
 \vdots & \vdots & \ddots & \vdots \\
1 & t_m & \ldots & t_m^m
\end{pmatrix}
\begin{pmatrix}
a_0\\
\vdots\\
a_m
\end{pmatrix}
=
\begin{pmatrix}
f(t_0)\\
\vdots\\
f(t_m)
\end{pmatrix}.
\end{equation*}
This system has a unique solution and so there exists a unique polynomial $g(X)$ of degree less than $m$ in $\R[X]$ such that $g(t_i)=f(t_i)$, for all $i$. From what we have just seen, there exists an element $u\in (t_0,t_m)$ such that 
$$
f^{(m)}(u)=g^{(m)}(u)=m!a_m \Longrightarrow a_m = \frac{f^{(m)}(u)}{m!}.
$$
Using Cramer's rule, we obtain
$$
a_m = \frac{W_m}{V_m} \Longrightarrow \frac{W_m}{V_m} = \frac{f^{(m)}(u)}{m!}.
$$
This ends the proof.\bx\\

The next preliminary result is interesting in that it shows that we often only need know that two functions have the same values at a limited number of points to establish that they have the same values at all points.

\begin{lemma}\label{lemHILBIRRED3}  Let $G$ be a connected open subset of $\C$ and $f$, $g$ analytic functions from $G$ into $\C$. Suppose that there exists a sequence $(z_n)_{n=1}^{\infty}$ in $G$and $z_0\in G$ such that
\begin{itemize}
\item for all $n\geq 1$, $z_n\neq z_0$;
\item $\lim _{n\rightarrow \infty} z_n=z_0$;
\item for all $n\geq 1$, $f(z_n)=g(z_n)$.
\end{itemize}
Then, for all $z\in G$, $f(z)=g(z)$.
\end{lemma}

\noindent \textsc{proof} As $f$ and $g$ are analytic on $G$, for all $z\in G$, we may find $r>0$ (depending on $z$) such that, for any $\tilde{z}$ with $\vert \tilde{z}-z\vert < r$, we have
$$
f(\tilde{z}) = \sum _{n=0}^{\infty}a_n(\tilde{z}-z^n)^n \qquad \mathrm{and}\qquad g(\tilde{z}) = \sum _{n=0}^{\infty}b_n(\tilde{z}-z^n)^n,
$$
where the coefficients $a_n$ and $b_n$ are in $\C$. Let $H$ be the subset of points $z\in G$ where the Taylor series of $f$ and $g$ coincide, i.e., $a_n=b_n$, for all $n\geq 0$.  We aim to show that $H=G$, which is sufficient to prove the result. First we show that $z_0\in H$ and so $H$ is nonempty. As $f$ and $g$ are analytic, we may find $r>0$ such that, for $\vert z-z_0\vert$, we have
$$
f(z) = \sum _{n=0}^{\infty}a_n(z-z_0)^n \quad \textrm{and} \quad \sum _{n=0}^{\infty}b_n(z-z_0)^n,
$$
where the coefficients $a_n$ and $b_n$ are in $\C$. We will use an induction argument to show that, for all $n$, $a_n=b_n$. To see this, we first notice that $f$ and $g$ are continuous at $z_0$, hence
$$
a_0 =f(z_0) = \lim _{n\rightarrow \infty}f(z_n) = \lim _{n\rightarrow \infty}g(z_n) =g(z_0) = b_0.
$$
Suppose now that we have established that $a_m=b_m$ up to a certain $m$. We set, for all $\vert z-z_0\vert \leq r$,  
$$
f^*(z) = a_{m+1} + a_{m+2}(z-z_0) + a_{m+3}(z-z_0)^2 + \cdots \quad \textrm{and} \quad g^*(z) = b_{m+1} + b_{m+2}(z-z_0) + b_{m+3}(z-z_0)^2 + \cdots.
$$
Then, for $0<\vert z-z_0\vert <r$, we have
$$
f^*(z) = \frac{1}{(z-z_0)^{m+1}}\left(f(z) -(a_0+a_1(z-z_0)+\cdots +a_m(z-z_0)^m)\right)
$$
and
$$
g^*(z) = \frac{1}{(z-z_0)^{m+1}}\left(g(z) -(b_0+b_1(z-z_0)+\cdots +b_m(z-z_0)^m)\right)
$$
and so, for all $n\geq 1$, $f^*(z_n)=g^*(z_n)$, since $z_n\neq z_0$ and $f(z_n)=g(z_n)$, for all $n\geq 1$. It follows that $a_{m+1}=b_{m+1}$, because $f^*$ and $g^*$ are continuous at $z$.  By induction, for all $n$, $a_n=b_n$. Hence the Taylor series of $f$ and $g$ coincide at $z_0$, which implies that $H\neq \emptyset$

The statement that $z\in H$ is equivalent to saying that $f^{(n)}(z)=g^{(n)}(z)$ for all $n\geq 0$. As the functions $f^{(n)}$ and $g^{(n)}$ are continuous, $H$ is closed in $G$. If $z\in H$, then there exists $r>0$ such that $f=g$ on the open disk $D(z,r)$. For any $\zeta \in D(z,r)$, there is a neighbourhood $N$ of $\zeta$ on which $f=g$ and so have the same Taylor series at $\zeta$. Therefore $D(z,r)\subset H$ and it follows that $H$ is open in $G$. Since $G$ is connected, we must have $H=G$ and so $f(z)=g(z)$, for all $z\in G$.\bx\\

In proving Hilbert's irreducibility theorem which follows (Theorem \ref{thmHILBIRRED1}), we will use two more preliminary results whose proofs we leave as exercises.

\begin{exercise}\label{exerHILBIRRED1} Suppose that $g\in \C[X]$, with $\deg g =m$, and $m+1$ integers $t_i$ such that $g(t_i)\in \Q$. Show that $g\in \Q[X]$.
\end{exercise}

\begin{exercise}\label{exerHILBIRRED2} Let $m\geq 1$ and $f(X)=\sum _{i=0}^ma_iX^i\in \Z[X]$. Show that, if $\alpha =\frac{p}{q}$, with $(p,q)=1$, is a root of $f$, then $p\vert a_0$ and $q\vert a_m$. Deduce that, if $f$ is monic, then any rational root of $f$ is an integer.
\end{exercise}

We now turn to the proof of Hilbert's result. 

\begin{theorem} (Hilbert's irreducibility Theorem)\label{thmHILBIRRED1} If $f(X,Y)\in \Q[X,Y]$ is irreducible, then there exists an infinite number of rational numbers $b$ such that $f_b(Y)=f(b,Y)$ is irreducible in $\Q[Y]$.
\end{theorem}

\noindent \textsc{proof} The proof is rather long and detailed, so we will proceed by steps.\\

\noindent \underline{1. The coefficient functions $y_j$}: From Proposition \ref{propHILBERTreduc1}, we know that all but a finite number number of values $b\in \Q$ are regular values of $f$. Consequently, we may choose $s_0\in \Q$, a regular value of $f$. Lemma \ref{lemHILBIRRED1} garantees the existence of $n$ roots $u_1(s), \ldots , u_n(s)$ of $f(s,Y)$, which are analytic functions on a $\C$-neighbourhood $N$ of $s_0$, which we may suppose to be connected. (This is of importance later.) 
As usual we write
$$
f(X,Y) = a_0(X) + a_1(X)Y + a_2(X)Y^2 + \cdots + a_n(X)Y^n,
$$
where the $a_i\in \Q[X]$, for $i=0,1,\ldots ,n$ , and $a_n(X)\neq 0$. Let us now consider 
$$
f(Y) = a_0 + a_1Y + a_2Y^2 + \cdots + a_nY^n,
$$
where, for $i=0,1,\ldots ,n$, $a_i$ is the polynomial function (with coefficients in $\Q$) associated with the polynomial $a_i(X)$. These functions are defined on $N$. Clearly, $f\in {\cal F}[Y]$, where ${\cal F}$ is the ring of polynomial functions defined on $N$, with coefficients in $\Q$. Moreover, it is clear that the functions $u_1, \ldots ,u_n$ are roots of $f$ and belong to the ring ${\cal A}$ of analytic functions defined on $N$. We have
$$
f(Y) = a_n\prod _{i=1}^n(-u_i +Y).
$$
We now let $S$ be a proper subset of $\N_n=\{1,\ldots ,n\}$, i.e., $S\neq \emptyset, \N_n$ and write
$$
\alpha (Y) = \prod _{i=1}^n(-u_i+Y) \qquad \beta (Y) =\prod_{i\in S} (-u_i+Y) \qquad \gamma (Y) = \prod _{i\notin S}(-u_i+Y).
$$
Then $\alpha , \beta ,\gamma \in {\cal A}[Y]$ and $\alpha =\beta \gamma$. If the coefficients of both $\beta$ and $\gamma$ are in ${\cal F}$, then the polynomial $f(X,Y)$ can be written as a product of polynomials of degree at least one in $\Q[X][Y]$. To see this, it is sufficient to write the equalities satisfied by the coefficients $a_0, a_1,\ldots ,a_n$ in the equality $f(Y)=a_n\beta \gamma$, where $a_n\neq 0$. Indeed, let
$$
\beta (Y) =\prod_{i=0}^kb_iY^i \qquad \textrm{and} \qquad \gamma (Y) = \prod _{j=0}^lc_jY^j,
$$
where $b_0,\ldots ,b_k, c_0,\ldots c_l\in {\cal F}$. Then
$$
a_0(s)= a_n(s) b_0(s)c_0(s),
$$
for an infinite number of rational numbers $s$, hence
$$
a_0(X) = a_n(X)b_0(X)c_0(X).
$$
Also,
$$
a_1(s) = a_n(s)\left(b_0(s)c_1(s)+b_1(s)c_0(s)\right),
$$
for an infinite number of rational numbers $s$, hence
$$
a_1(X) = a_n(X)\left(b_0(X)c_1(X) + b_1(X)c_0(X)\right).
$$
Continuing in the same way, we find that 
$$
a_k(X) = a_n(X)\sum _{i+j=k}b_i(X)c_j(X),
$$
for $k=0, 1, \ldots ,n$. If we set 
$$
\beta (X,Y) = \sum _{i=0}^kb_i(X)Y^i \qquad\textrm{and}\qquad \gamma (X,Y) = \sum _{j=0}^lc_j(X)Y^j,
$$
then 
$$
f (X,Y) = a_n(X)\beta (X,Y)\gamma (X,Y),
$$
which implies that $f(X,Y)$ is not irreducible, a contradiction. It follows that, for any proper subset $S$ of $\N_n$, either $\beta$ or $\gamma$ has a coefficient $y$ which is not a polynomial function with rational coefficients. 

If we replace the $b_i$ and $c_j$ with quotients of polynomial functions with rational coefficients, then analogous calculations to those which we have just employed show that $f(X,Y)$ is reducible in $F(X)[Y]$, which from Gauss's Lemma is not possible, because $f(X,Y)$ is irreducible in $F[X][Y]$. Therefore we may assume that $y$ is not even a quotient of polynomial functions with rational coefficients. We number the distinct functions functions $y_1, \ldots ,y_{2^n-2}$. (We do not say that these functions are distinct; certain of them may be the same.)\\

\noindent \underline{2. A condition for the irreducibility of $f(s,Y)$}: Suppose that $s\in N\cap \Q$ and that $y_1(s), \ldots ,y_{2^n-2}(s)$ are all in $\C\setminus \Q$. We claim that $f(s,Y)$ is irreducible in $\Q[Y]$. Indeed, we can always write
$$
f(s,Y) =a_n(s) \prod _{i\in S}\left(-u_i(s) +Y\right)\prod_{i\notin S}\left(-u_i(s) +Y\right),
$$
for any proper subset $S$ of $\N_n$. As $s$ is rational, $a_n(s)$ is also rational. If we evaluate the coefficients of $\beta$ and $\gamma$ at $s$, we obtain the coefficients of $\prod _{i\in S}\left(-u_i(s) +Y\right)$ and $\prod_{i\notin S}\left(-u_i(s) +Y\right)$. By the choice of $s$, at least one of these coefficients is not rational. Thus, $f(s,Y)$ is irreducible in $\Q[Y]$. In order to prove our theorem, it is sufficient to obtain an infinite number of elements  $s\in N\cap \Q$ such that $y_1(s), \ldots ,y_{2^n-2}(s)$ are all in $\C\setminus \Q$.\\

\noindent \underline{3. Studying the functions $y_j$}: We aim to look at the functions $y_j$ in more detail. To simplify the notation, let us write $y$ for $y_j$. We notice that there exists $T'>0$ such that $s_0+\frac{1}{t}\in N$, if $t> T'$. We define the function $\delta$ on $(T',\infty )$ by
$$
\delta (t) = y(s_0+\frac{1}{t}).
$$
Let us denote ${\cal G}$ the set of functions defined on $(T',\infty )$ by
$$
\xi (t) = \frac{h(t)}{g(t)},
$$
where $h$ and $g$ are polynomial functions with rational coefficients and $g$ is not the zero function. Clearly, ${\cal G}$ is a field. We claim that $\delta$ is algebraic over ${\cal G}$. To see this, first let us define the function $v_i$ on $(T',\infty )$ by
$$
v_i(t) = u_i (s_0 + \frac{1}{t}).
$$
Then, for all $t>T'$, we have
$$
a_0(s_0+\frac{1}{t}) + a_1(s_0+\frac{1}{t})v_i(t) + a_2(s_0+\frac{1}{t})v_i(t)^2 + \cdots + a_n(s_0+\frac{1}{t})v_i(t)^n = 0.
$$
Multiplying by the appropriate power of $t$, we obtain the expression
$$
h_0(t) + h_1(t)v_i(t) + h_2(t)v_i(t)^2 + \cdots + h_n(t)v_i(t)^n,
$$
where $h_0, \ldots ,h_n$ are polynomial functions defined on $(T',\infty )$, with coefficients in $\Q$. Hence,
$$
h_0 + h_1v_i + h_2v_i^2 + \cdots + h_nv_i^n
$$
is the zero function and so $v_i$ is algebraic over ${\cal G}$. As the algebraic elements over a field form a field, $\delta$ is algebraic over ${\cal G}$ and so is the root of an equation
$$
d_0 + d_1H \ldots + \ldots d_mH^m= 0,
$$
where $d_0, \ldots ,d_m\in {\cal G}$ and $0$ denotes the zero function.  We may suppose that the $d_i$ are polynomial functions. (It is sufficient to multiply by the product of the denominators of the $d_i$, if necessary.) We may even suppose that the coefficients of the polynomial functions $d_i$ are integers. We now multiply the coefficients of the equation satisfied by $\delta$ by $d_m^{m-1}$ to obtain
$$
d_0d_m^{m-1} + d_1d_m^{m-2}d_m\delta + d_2d_m^{m-3}(d_m\delta )^2 +\cdots + d_{m-1}(d_m\delta )^{m-1} + (d_m\delta )^m = 0.
$$
Thus $z=d_m\delta$ is a root of the polynomial
\begin{equation}\label{eqnHILBIRRED2}
g(Z) = b_0 + b_1Z + \cdots + b_{m-1}Z^{m-1} + Z^m,
\end{equation}
where $b_j=d_jd_m^{m-j-1}$, for $j=0, 1, \ldots ,m-1$. Clearly, the coefficients of $g(Z)$ are integers. We claim that, if $t\in \Z$, with $t>T'$, and $\delta (t)\in \Q$, then $z(t)\in \Z$. Indeed, $z(t)=d_m(t)\delta (t)$ implies that $z(t)\in \Q$. Also, $z(t)$ is a root of the equation
$$
b_0(t) + b_1(t)Z + \cdots + b_{m-1}(t)Z^{m-1} + Z^m,
$$
which is a monic polynomial with coefficients in $\Z$. From Exercise \ref{exerHILBIRRED2}, $z(t)$ is an integer.\\

\noindent \underline{4. Studying the functions $z_j$}: To simplify the notation, we will write $z$ for $z_j$: Our next step is to show that there are relatively few integers $t>T'$ such that $z(t)$ is an integer. If this is the case, then we may find many integers $t$ such that $z(t)$ is not an integer. For such $t$, $\delta (t)$ cannot be rational, which implies that $y(s_0+\frac{1}{t})$ is not rational.

Lemma \ref{lemHILBIRRED1} ensures us that, for $i=1,\ldots ,n$, $u_i$ is an analytic function on the neighbourhood $N$ of $s_0$. As sums and products of analytic functions are analytic, for $j=1,\ldots ,2^n-2$, $y_j$ is analytic on $N$. Reducing the size of $N$ to a neighbourhood $N'$ of $s_0$ if necessary, for $s_0+x\in N'$, we may write 
$$
y(s_0+x)= e_0 + e_1x + e_2x^2 + \ldots + e_kx^k + \ldots, 
$$
where the coefficients $e_i\in \C$. There exists $T''\geq T'$ such that, if $t>T''$, then $y(s_0+\frac{1}{t})\in N'$ and so
$$
y(s_0+\frac{1}{t}) = e_0 + e_1\frac{1}{t} + e_2\left(\frac{1}{t}\right)^2 + \ldots + e_k\left(\frac{1}{t}\right)^k + \ldots.
$$
As $d_m$ is a polynomial, we may write
$$
z(t) = d_m(t)\delta (t) = d_m(t)y(s_0+\frac{1}{t}) = c_lt^l + \cdots + c_1t + c_0 + c_{-1}t^{-1} + \cdots + c_{-k}t^{-k} +\cdots ,
$$
with $c_i\in \C$.  

There are three possibilities:
\begin{itemize}
\item {\bf a.} $z$ is a polynomial function;
\item {\bf b.} $z$ is not a polynomial function and has at least one coefficient $c_i\in \C\setminus \R$;
\item {\bf c.} $z$ is not a polynomial function and all the coefficients of $z$ are real.
\end{itemize}

We consider the first case. We claim that at least one of the coefficients must be in $\C\setminus \Q$. If this is not the case, then $\delta (t)=\frac{z(t)}{d_m(t)}$, for $t>T''$. Let $(t_n)$ be sequence of values of $t>T''$ converging to $\infty$. If we set $s_n = s_0+\frac{1}{t_n}$, then the numbers $s_n$ converge to $s_0$ and
$$
y(s_n) = y(s_0+\frac{1}{t_n}) =\delta (t_n) = \frac{z(t_n)}{d_m(t_n)} = \frac{z((s_n-s_0)^{-1})}{d_m((s_n-s_0)^{-1})}.
$$
If we multiply both $z$ and $d_m$ by an appropriate power of $s_n-s_0$, then we may find polynomial functions with rational coefficients $\hat{z}$ and $\hat{d}_m$ such that 
$$
\frac{z((s_n-s_0)^{-1})}{d_m((s_n-s_0)^{-1})}= \frac{\hat{z}(s_n)}{\hat{d}_m(s_n)}.
$$
From Lemma \ref{lemHILBIRRED3}, we obtain that $y=\frac{\hat{z}}{\hat{d}_m}$, a contradiction. (We can apply Lemma \ref{lemHILBIRRED3}, because $y$ and $\frac{\hat{z}}{\hat{d}_m}$ are defined on the connected set $N$.) This proves our claim. Thus at least one coefficient of $z$ belongs to the set $\C\setminus \Q$. Consequently, from Exercise \ref{exerHILBIRRED1}, there can only be a finite number of integers $t$ such that $z(t)$ is an integer. In this case we may choose $T'''>T''$ such that $z(t)$ is not an integer, if $t>T'''$.\\

Let us now consider the second case. Suppose that $i_0$ is the largest subscript $i$ for which $c_i\in \C\setminus \R$. Then
$$
\lim _{t\rightarrow \infty}\frac{\Img z(t)}{t^{i_{0}}} = \Img c'_{i_0}\neq 0.
$$ 
Hence, we may find $T'''\geq T''$ such that $\frac{z(t)}{t^{i_0}}\notin \R$, for $t>T'''$. This implies that $z(t)\notin \R$ and so is not an integer, for $t>T'''$. \\

The third case is more difficult to handle. Here all the coefficients $c_i$ are real and at least one coefficient $c_i$, with $i$ negative, is nonzero. By differentiating $z$ a sufficient number of times we can eliminate all nonnegative powers of $t$ to obtain
$$
z^{(m)}(t) = pt ^{-q}+\cdots , 
$$
where $p$ is a nonzero real number, $q$ a positive integer greater than $m$ and the dots represent terms of higher powers of $t^{-1}$. As
$$
\lim _{t\rightarrow \infty}t^qz^{(m)}(t) = p,
$$
there exists $T'''\geq T''$ such that
$$
t> T''' \Longrightarrow 0<\vert z^{(m)}(t)\vert \leq 2\vert p\vert t^{-q}.
$$
Now we use Lemma \ref{lemHILBIRRED2}. Let $t_0<t_1<\cdots <t_m$ be integers such that $t_i\geq T'''$ and $z(t_i)\in \Z$, for all $i$. For a certain number $u\in (t_0,t_m)$ we have
$$
\frac{2\vert p\vert}{m!t_0^q} > \frac{2\vert p\vert}{m!u^q}\geq \frac{\vert z^{(m)}(u)\vert}{m!} = \frac{\vert W_m\vert}{\vert V_m\vert }.
$$
As $z^{(m)}(u)\neq 0$, $W_m\neq 0$, which implies that $W_m$ is a positive integer and so $\vert W_m\vert \geq 1$. Therefore
$$
\frac{m!}{2\vert p\vert }t_0 ^q < \vert V_m\vert =\prod _{i>j}(t_i -t_j)<(t_m-t_0)^{\frac{m(m-1)}{2}}.
$$
This implies that there are positive constants $\alpha$ and $\beta$ such that $\alpha t_0^{\beta}<t_m-t_0$.

Now let $r$ be the number of distinct functions $z_j$ in this third case. Without loss of generality, let us suppose that these are the functions $z_1, \ldots ,z_r$. For each $j$, we have $m_j$, $\alpha _j$ and $\beta _j$ such that, if we have integers $t_0<t_1<\ldots <t_{m_j}$, with $z_j(t_i)\in \Z$, for $i=0, 1, \ldots ,m_j$, then $\alpha _jt_0^{\beta _j}<t_{m_j}-t_0$. We set $\bar{m}=\max m_j$ and take $U\in \Z$ such that $\alpha _jU^{\beta _j}\geq rm$. We now consider the interval $I=[U,U+rm]$. If $t_0<t_1<\ldots <t_{m_j}$ is a sequence of $m_j+1$ integers in $I$, then
 $$
\alpha _jt_0^{\beta _j}\geq \alpha _jU^{\beta _j}\geq rm = (t_0+rm)-t_0.
$$
This implies that $I$ contains at most $m_j$ integers $t$ such that $z_j(t)\in \Z$. 

If we now consider all the $z_j$ in the third case, we see that the interval contains at most $m_1+\ldots +m_r$ integers $t$ such that $z_j(t)\in \Z$ for some $j=1, \ldots ,r$, i.e., at most $rm$ integers $t$ such that $z_j(t)\in \Z$, for some $j=1,\ldots ,r$. However, $I$ contains $nm+1$ integers, so there is at least one integer $t\in I$ such that $z_j(t)\notin \Z$, for $j=1,\ldots ,r$. We may find an infinite number of such $t$ by taking a sequence of intervals $I_k=[U_k,U_k+rm]$, with $U_{k+1}>U_k+rm$.\\

\noindent \underline{5. The final step}: Using ou previous work, we show that there is an infinite sequence of integers $t$ such that $z_j(t)$ is not an integer, for all $z_j$. As we have seen, there are three possibilities for $z_j$. For those which fall in the categories {\bf a} or {\bf b.}, there is a number $T'''$ such that, if $t>T'''$, then $z_j(t)\notin \Z$. If we take $T''''$ equal to be the maximum of all such $T'''$, then $z_j(t)\notin \Z$, for those $z_j(t)$, where $z_j$ is in either category {\bf a.} or {\bf b.} If the category {\bf c.} is empty, then we have finished. 

If this is not the case and $z_1, \ldots ,z_r$ belong to the third case, then we can find a sequence of integers $t$ such that $z_j(t)\notin \Z$, for $j=1,\ldots ,r$. We may take these integers greater than $T''''$ and so we have an infinite sequence of integers $t$ such that $z_j(t)\notin \Z$ for all $j$. This finishes the proof.\bx\\

\noindent {\bf Remark} The rational numbers $b$, which we have found, such that $f_b(Y)$ is irreducible are of the form $s_0+t^{-1}$, where $t$ is a positive integer. Of course, there are certainly others: we only need to take $s_0'$ sufficiently far from $s_0$. 

\begin{exercise} In section $4$ of the above proof, we assumed that $W_m$ is an integer. Why is this so?
\end{exercise}

\noindent {\bf Remark} In this article we have been concerned with irreducible specializations. We have not considered reducible specializations. For many number fields $K$ the number of reducible specializations of irreducible polynomials in $K[X,Y]$ is infinite. However, there are number fields with irreducible polynomials for which this is not the case. For a recent discussion of this question see \cite{muller}.

\vspace{2cm}

\noindent \underline{\bf Basic results from Galois theory}

\begin{result}\label{ufd} Let $R$ be a unique factorization domain, with quotient field $F$, and $f\in R[X]$. Then, if $f$ is nonconstant and irreducible in $R[X]$, then $f$ is irreducible in $F[X]$. On the other hand, if $f$ is primitive and irreducible in $F[X]$, then $f$ is irreducible in $R[X]$.
\end{result} 

\begin{result}\label{lemSPLIT1} Let $f\in F[X]$ be irreducible and $E$ an extension of $F$ which contains a root $\alpha$ of $f$. Then there is an isomorphism
$$
\Phi : F[X]/(f) \longrightarrow F(\alpha )
$$
which fixes $F$, i.e., for $g$ constant, $\Phi (g+(f))=g$, and such that $\Phi (X+(f))=\alpha$.
\end{result}

\begin{result}\label{thSPLIT2} Let $F$ and $F'$ be fields, $\sigma :F\longrightarrow F'$ an isomorphism, $f\in F[X]$ and $f^*\in F'[X]$ the polynomial corresponding to $f$. If $E$ is a splitting field of $f$ and $E'$ a splitting field of $f^*$, then there is an isomorphism $\tilde{\sigma}:E\longrightarrow E'$ extending $\sigma$.
\end{result}

\begin{result}\label{propSPLIT2} Let $\sigma :F\longrightarrow F'$ be an isomorphism and $f\in F[X]$ irreducible. If $E$ (resp. $E'$) is an extension of $F$ (resp. $F'$) and $\alpha $ (resp. $\alpha '$) a root of $f$ (resp. $f^*$) in $E$ (resp. $E'$), then there is an isomorphism $\hat{\sigma} : F(\alpha )\longrightarrow F'(\alpha ')$ extending $\sigma$, with $\hat{\sigma}(\alpha )=\alpha '$. This isomorphism is unique.
\end{result}

\begin{result}\label{thSEPext4} Let $E$ be a finite separable extension of a field $F$ of degree $n$. Then the field of fractions $E(X)$ is a finite extension of degree $n$ of the field of fractions $F(X)$.
\end{result}

\begin{result}\label{NORMALprop2} Suppose that $K/F$ and $E/K$, with $E$ normal over $F$. Then $E$ is normal over $K$.
\end{result}
 
\begin{result}\label{NORMALth1} The finite extension $E$ of $F$ is normal if and only if $E$ is the splitting field of a polynomial $f\in F[X]$.
\end{result}

\begin{result}\label{NormalClos1} If $E$ is  finite extension of the field $F$ and $N$ the normal closure of $E$ over $F$, then $N$ is a finite extension of $F$.
\end{result}

\begin{result}\label{thGALGRP1} If $E$ is a finite Galois extension of $F$, then we have $\vert Gal(E/F)\vert =[E:F]$.
\end{result}

\begin{result}\label{thGALPOLYirred1} Let $f$ be a separable polynomial in $F[X]$ of degree $n$ with Galois group $G=Gal(E/F)$. If $f\in F[X]$ is irreducible, then the action of $G$ on the set of roots of $f$ is transitive.
\end{result}

\end{document}